# On a multiplicative order of Gauss periods and related questions

Roman Popovych

December 19, 2010

**Abstract.** We obtain explicit lower bounds on multiplicative order of elements that have more general form than finite field Gauss period. In a partial case of Gauss period this bound improves the previous bound of O.Ahmadi, I.E.Shparlinski and J.F.Voloch

## 1    Introduction

A problem of constructing primitive element for the given finite field is notoriously difficult in computational number theory. That is why one considers less restrictive question: to find an element with large order. It is sufficient in this case to obtain a lower bound of the order.

Such question in particular arises in concern with primality proving AKS algorithm. This unconditional, deterministic polynomial-time algorithm for proving primality of integer $n$ was presented by M.Agrawal, N.Kayal and N.Saxena in [1].

The idea of AKS consists in the following: to choose small integer $r$ with specific properties and to show that the set of elements $\theta+a$, $a=1,\ldots,\lfloor\sqrt{r}\log n\rfloor$ generates large enough subgroup in the multiplicative group of the finite field $Z_p[x]/h(x)$ where $\theta$ is a coset of $x$, $p$ is a prime divisor of $n$, $h(x)$ is an irreducible divisor of $r$-th cyclotomic polynomial $\Phi_r(x)=C_r(x)=x^{r-1}+x^{r-2}+\ldots+x+1$. Usually in AKS implementation $r$ is prime and $h(x)=\Phi_r(x)$.

The fastest known deterministic version of AKS runs in $(\log n)^{6+o(1)}$ time, the fastest randomized version – in $(\log n)^{4+o(1)}$ time [6]. If given in [1, p.791] conjecture were true, this would improve a complexity of AKS to $O(\log n)^{3+o(1)}$. The conjecture means that the element $\theta-1$ generates large enough subgroup.

We denote $F_q$ finite field with $q$ elements where $q$ is a power of prime number $p$. Let $r=2s+1$ be a prime number coprime with $q$. Let $q$ be a primitive root modulo $r$, that is a multiplicative order of $q$ modulo $r$ equals to $r-1$.

Set $F_q(\theta)=F_{q^{r-1}}=F_q[x]/\Phi_r(x)$ where $\theta=x(\mathrm{mod}\,\Phi_r(x))$, $\theta$ generates an extension $F_q(\theta)=F_{q^{r-1}}$. It is clear that $\theta^r=1$. The element $\theta+\theta^{-1}$ is called a Gauss period of type ((r-1)/2,2). It allows to construct normal base [2].

$<u_1,\ldots,u_k>$ denotes a group generated by elements $u_1,\ldots,u_k$. $\rho_2(n)$ denotes the highest power of 2 dividing integer $n$.

A partition of integer C consists of such non-negative integers $u_1,\ldots,u_C$ that $\sum_{j=1}^{C} ju_j = C$. $U(C)$ denotes a number of partitions of $C$. $U(C,d)$ denotes a number of such partitions $u_1,\ldots,u_C$ of $C$ for which $u_1,\ldots,u_C \leq d$. $Q(C,d)$ denotes a number of such partitions $u_1,\ldots,u_C$ of $C$ for which $u_1,\ldots,u_C \not\equiv 0 \bmod d$.

It is shown in [2] that a multiplicative order of Gauss period $\theta+\theta^{-1}$ which generates normal base over finite field is at least $U((r-3)/2,p-1)$.

We generalize this result and show that for any integer $e$, any integer $f$ coprime with $r$, any non-zero element $a$ in the field $F_q$ the multiplicative order of the element $\theta^e(\theta^f+a)$ in the field $F_q(\theta)$ is at least $U(r-2,p-1)$. In particular, a multiplicative order of the element $\theta+\theta^{-1}=\theta^{-1}(\theta^2+1)$ is at least $U(r-2,p-1)$. This bound improves the previous bound of O.Ahmadi, I.E.Shparlinski and J.F.Voloch given in [2]. We show that if $a \neq 1,-1$ then a multiplicative order of element $\theta^e(\theta^f+a)$ is least $[U((r-3)/2,p-1)]^2/2$. We also prove that an order of subgroup $<\theta+\theta^{-1}, (a\theta+1)(\theta+a)^{-1}>$ is at least $[U(r-2,p-1)\ U((r-3)/2,p-1)]/2$ and construct a generator of this cyclic subgroup.

We also give explicit lower bounds for the multiplicative order of the elements using results from [3,7].

A construction of large order elements in the case $q \equiv 1 \bmod r$ ($r$ is not primitive modulo $r$) is given in [5]. Overview of some alternative constructions of large order elements in a finite field can be found in [4].

## 2      Multiplicative orders of finite field elements

**Theorem 1.** *Let $q$ be a power of prime number $p$, $r=2s+1$ be a prime number coprime with $q$, $q$ be a primitive root modulo $r$, $\theta$ generates an extension $F_q(\theta)=F_{q^{r-1}}$, $e$ be any integer, $f$ be any integer coprime with $r$, $a$ be any non-zero element in the finite field $F_q$. Then*

*(a) element $\theta^e(\theta^f+a)$ has a multiplicative order at least $U(r-2,p-1)$*

*(b) if $a \neq 1,-1$ then element $\theta^e(\theta^f+a)$ has a multiplicative order at least $[U((r-3)/2,p-1)]^2/2$*

*Proof.* (*a*) Let us consider automorphism of the field $F_q(\theta)$ which takes $\theta$ to $\theta^f$. Since the automorphism sends element $\theta^g(\theta+a)$ to element $\theta^e(\theta^f+a)$, where $g \equiv ef^{-1} \bmod r$, multiplicative orders of these elements coincide. Hence, it is sufficient to prove that element of the form $\theta^e(\theta+a)$ has a multiplicative order at least $U(r-2,p-1)$.

Since $q$ is primitive modulo $r$, then for each $j=1,\ldots,r-2$ exists such integer $\alpha_j$ that $q^{\alpha_j} \equiv j \bmod r$. Then $q^{\alpha_j}$-th powers of element $\theta^e(\theta+a)$ which are equal to

$$\left(\theta^e(\theta+a)\right)^{q^{\alpha_j}} = \theta^{eq^{\alpha_j}}(\theta^{q^{\alpha_j}}+a) = \theta^{ej}(\theta^j+a)$$

belong to group $<\theta^e(\theta+a)>$. We consider the following products

$$\prod_{j=1}^{r-2}[\theta^{ej}(\theta^j+a)]^{u_j}, \text{ where } \sum_{j=1}^{r-2} ju_j = r-2,\ 0 \le u_1,...,u_{r-2} \le p-1,$$

which also belong to the group, and show that a number of such products equals to $U(r-2, p-1)$.

Let partitions $(u_1,...,u_{r-2})$ and $(v_1,...,v_{r-2})$ of integer $r$-2, where it part appears no more than $p$-1 times, be different, but the correspondent products coincide:

$$\prod_{j=1}^{r-2}[\theta^{ej}(\theta^j+a)]^{u_j} = \prod_{j=1}^{r-2}[\theta^{ej}(\theta^j+a)]^{v_j} \tag{1}$$

Then we have from (1):

$$\theta^{e\sum_{j=1}^{r-2} ju_j}\prod_{j=1}^{r-2}(\theta^j+a)^{u_j} = \theta^{e\sum_{j=1}^{r-2} jv_j}\prod_{j=1}^{r-2}(\theta^j+a)^{v_j}$$

$$\prod_{j=1}^{r-2}(\theta^j+a)^{u_j} = \prod_{j=1}^{r-2}(\theta^j+a)^{v_j} \tag{2}$$

Since the characteristic polynomial of $\theta$ is the polynomial $\Phi_r(x)$ we obtain

$$\prod_{j=1}^{r-2}(x^j+a)^{u_j} = \prod_{j=1}^{r-2}(x^j+a)^{v_j} \tag{3}$$

As there are polynomials of degree $r$-2<deg $\Phi_r(x)$ from left and right side in equality (3) then these polynomials coincide as polynomials over $F_q$.

Let $k$ be the smallest integer for which $u_k \ne v_k$. Without loss of generality suppose $u_k > v_k$. After removing from (3) common factors we obtain

$$(x^k+a)^{u_k-v_k}\prod_{j=k+1}^{r-2}(x^j+a)^{u_j} = \prod_{j=k+1}^{r-2}(x^j+a)^{v_j} \tag{4}$$

Then there is the term $(u_k-v_k)a^{u_k-v_k-1}a^{\sum_{j=k+1}^{r-2} u_j}x^k$ in the polynomial in the left side of (4) with minimal non-trivial power of $x$. Since $0<u_k,v_k<p$-1, $u_k \ne v_k$, $a \ne 0$, the term is non-zero. And all terms with non-trivial power of $x$ in the polynomial in the right side have a power higher than $k$ - a contradiction.

(*b*) Order of a multiplicative group of the field $F_{q^{r-1}}$ equals to $q^{r-1}$-1=$(q^{(r-1)/2}$-1)$(q^{(r-1)/2}$+1). The factors $q^{(r-1)/2}$-1 and $q^{(r-1)/2}$+1 have the greatest common divisor 2 since their sum equals to $2q^{(r-1)/2}$.

Subgroup $H$ of the group generated by element $\theta^e(\theta+a)$ contains subgroup $H_1$ generated by element

$$v = [\theta^e(\theta+a)]^{q^{(r-1)/2}-1} = \theta^{eq^{(r-1)/2}}(\theta^{q^{(r-1)/2}}+a)(\theta+a)^{-1} = \theta^{-e}(\theta^{-1}+a)(\theta+a)^{-1} = \theta^{-(e+1)}(a\theta+1)(\theta+a)^{-1}$$

and subgroup $H_2$ generated by element

$w = [\theta^e(\theta+a)]^{q^{(r-1)/2}+1} = \theta^{eq^{(r-1)/2}}(\theta^{q^{(r-1)/2}}+a)(\theta+a) = \theta^{-e}(\theta^{-1}+a)(\theta+a) = \theta^{-(e+1)}(a\theta+1)(\theta+a)$

(note that since $q$ is primitive modulo $r$ and $r$ is prime, then $q^{r-1}=1$ mod $r$, a $q^{(r-1)/2}=-1$ mod $r$).

Element $[\theta^e(\theta+a)]^{q^{(r-1)/2}-1}$ has order that divides $q^{(r-1)/2}+1$, and element $[\theta^e(\theta+a)]^{q^{(r-1)/2}+1}$ has order that divides $q^{(r-1)/2}-1$.

Let us consider element $z = \begin{cases} v^2 w \text{ if } \rho_2(q^{(r-1)/2}-1) = 2 \\ vw^2 \text{ if } \rho_2(q^{(r-1)/2}+1) = 2 \end{cases}$

If $\rho_2(q^{(r-1)/2}-1) = 2$ then orders of elements $v^2$ and $w$ are coprime. If $\rho_2(q^{(r-1)/2}+1) = 2$ then orders of elements $v$ and $w^2$ are coprime. In both cases an order of element $z$ is a product of orders of $v$ and $w$ divided by 2.

Element $w = \theta^{-(e+1)}(a\theta+1)(\theta+a)$ generates subgroup $H_2$ of order at least $U((r-3)/2, p-1)$.

Indeed, elements $\theta^{-2(e+1)}(a\theta^2+1)(\theta^2+a), \ldots, \theta^{-[(r-3)/2](e+1)}(a\theta^{(r-3)/2}+1)(\theta^{(r-3)/2}+a)$ belong to subgroup $H_2$.

We consider products of the powers

$$\prod_{j=1}^{(r-3)/2} [\theta^{-j(e+1)}(a\theta^j+1)(\theta^j+a)]^{u_j}, \text{ where } \sum_{j=1}^{(r-3)/2} ju_j = (r-3)/2, 0 \le u_1, \ldots, u_{r-2} \le p-1,$$

which also belong to the subgroup, and show that a number of such products equals to $U((r-3)/2, p-1)$.

Let partitions $(u_1, \ldots, u_{r-2})$ and $(v_1, \ldots, v_{r-2})$ of integer $(r-3)/2$, where it part appears no more than $p-1$ times, be different, but the correspondent products coincide:

$$\prod_{j=1}^{(r-3)/2} [\theta^{-j(e+1)}(a\theta^j+1)(\theta^j+a)]^{u_j} = \prod_{j=1}^{(r-3)/2} [\theta^{-j(e+1)}(a\theta^j+1)(\theta^j+a)]^{v_j} \qquad (5)$$

Then we have from (5):

$$\theta^{-(e+1)\sum_{j=1}^{(r-3)/2} ju_j} \prod_{j=1}^{(r-3)/2} [(a\theta^j+1)(\theta^j+a)]^{u_j} = \theta^{-(e+1)\sum_{j=1}^{(r-3)/2} jv_j} \prod_{j=1}^{(r-3)/2} [(a\theta^j+1)(\theta^j+a)]^{v_j}$$

$$\prod_{j=1}^{(r-3)/2} [(a\theta^j+1)(\theta^j+a)]^{u_j} = \prod_{j=1}^{(r-3)/2} [(a\theta^j+1)(\theta^j+a)]^{v_j} \qquad (6)$$

Then we can write

$$\prod_{j=1}^{(r-3)/2} [(ax^j+1)(x^j+a)]^{u_j} = \prod_{j=1}^{(r-3)/2} [(ax^j+1)(x^j+a)]^{v_j} \qquad (7)$$

There are polynomials of degree $r-3 < \deg \Phi_r(x)$ from left and right side in equality (7).

Let $k$ be the smallest integer for which $u_k \neq v_k$. Without loss of generality suppose $u_k > v_k$. After removing from (7) common factors we obtain

$$[(ax^k+1)(x^k+a)]^{u_k-v_k} \prod_{j=k+1}^{(r-3)/2}[(ax^j+1)(x^j+a)]^{u_j} = \prod_{j=k+1}^{(r-3)/2}[(ax^j+1)(x^j+a)]^{v_j}$$

$$[ax^{2k}+(a^2+1)x^k+a]^{u_k-v_k} \prod_{j=k+1}^{(r-3)/2}[(ax^j+1)(x^j+a)]^{u_j} = \prod_{j=k+1}^{(r-3)/2}[(ax^j+1)(x^j+a)]^{v_j} \quad (8)$$

Then there is the term $(u_k-v_k)(a^2+1)a^{u_k-v_k-1}a^{\sum_{j=k+1}^{(r-3)/2} u_j} x^k$ in the polynomial in the left side of (8) with minimal non-trivial power of $x$. Since $0 < u_k, v_k < p-1$, $u_k \neq v_k$, $a \neq 0, 1, -1$, the term is non-zero. And all terms with non-trivial power of $x$ in the polynomial in the right side have a power higher than $k$ - a contradiction.

Element $v = \theta^{-(e+1)}(a\theta+1)(\theta+a)^{-1}$ generates subgroup $H_1$ of order at least $U((r-3)/2, p-1)$.

Indeed, elements $\theta^{-2(e+1)}(a\theta^2+1)(\theta^2+a)^{-1}, \ldots, \theta^{-[(r-3)/2](e+1)}(a\theta^{(r-3)/2}+1)(\theta^{(r-3)/2}+a)^{-1}$ belong to subgroup $H_1$.

We consider products of the powers

$$\prod_{j=1}^{(r-3)/2}[\theta^{-j(e+1)}(a\theta^j+1)(\theta^j+a)^{-1}]^{u_j}, \text{ where } \sum_{j=1}^{(r-3)/2} ju_j = (r-3)/2, \ 0 \le u_1, \ldots, u_{r-2} \le p-1,$$

which also belong to the subgroup, and show that a number of such products equals to $U((r-3)/2, p-1)$.

Let partitions $(u_1, \ldots, u_{r-2})$ and $(v_1, \ldots, v_{r-2})$ of integer $(r-3)/2$, where it part appears no more than $p$-1 times, be different, but the correspondent products coincide:

$$\prod_{j=1}^{(r-3)/2}[\theta^{-j(e+1)}(a\theta^j+1)(\theta^j+a)^{-1}]^{u_j} = \prod_{j=1}^{(r-3)/2}[\theta^{-j(e+1)}(a\theta^j+1)(\theta^j+a)^{-1}]^{v_j} \quad (9)$$

Then we have from (9):

$$\theta^{-(e+1)\sum_{j=1}^{(r-3)/2} ju_j} \prod_{j=1}^{(r-3)/2}[(a\theta^j+1)(\theta^j+a)^{-1}]^{u_j} = \theta^{-(e+1)\sum_{j=1}^{(r-3)/2} jv_j} \prod_{j=1}^{(r-3)/2}[(a\theta^j+1)(\theta^j+a)^{-1}]^{v_j}$$

$$\prod_{j=1}^{(r-3)/2}[(a\theta^j+1)(\theta^j+a)^{-1}]^{u_j} = \prod_{j=1}^{(r-3)/2}[(a\theta^j+1)(\theta^j+a)^{-1}]^{v_j} \quad (10)$$

Then we obtain

$$\prod_{j=1}^{(r-3)/2}[(ax^j+1)(x^j+a)^{-1}]^{u_j} = \prod_{j=1}^{(r-3)/2}[(ax^j+1)(x^j+a)^{-1}]^{v_j} \quad (11)$$

There are polynomials of degree $r-3 < \deg \Phi_r(x)$ from left and right side in equality (11).

Let $k$ be the smallest integer for which $u_k \neq v_k$. Without loss of generality suppose $u_k > v_k$. After removing from (11) common factors we obtain

$$[(ax^k+1)(x^k+a)^{-1}]^{u_k-v_k} \prod_{j=k+1}^{(r-3)/2} [(ax^j+1)(x^j+a)^{-1}]^{u_j} = \prod_{j=k+1}^{(r-3)/2} [(ax^j+1)(x^j+a)^{-1}]^{v_j}$$

$$(ax^k+1)^{u_k-v_k} \prod_{j=k+1}^{(r-3)/2} (ax^j+1)^{u_j}(x^j+a)^{v_j} = (x^k+a)^{u_k-v_k} \prod_{j=k+1}^{(r-3)/2} (ax^j+1)^{v_j}(x^j+a)^{u_j} \qquad (12)$$

Let us denote absolute term for $\prod_{j=k+1}^{(r-3)/2}(ax^j+1)^{u_j}(x^j+a)^{v_j}$ by $c$, and absolute term for

$\prod_{j=k+1}^{(r-3)/2}(ax^j+1)^{v_j}(x^j+a)^{u_j}$ by $d$. It is clear that $c \neq 0$, $d \neq 0$. Since absolute terms from left and from right in (12) must be equal, we have $c = a^{u_k-v_k}d$. Since coefficients near $x^k$ from left and from right in (12) must be equal, we have $c(u_k-v_k)a = d(u_k-v_k)a^{u_k-v_k-1}$. Then $a^2 = 1$ - a contradiction.

Hence, order of element $\theta^e(\theta+a)$ is at least $[U((r-3)/2, p-1)]^2/2$. □

**Remark 1.** Element $\theta + \theta^{-1}$ belongs to the subfield $F_{q^{(r-1)/2}}$ of the field $F_q(\theta) = F_{q^{r-1}}$.

**Corollary 2.** *Element $\theta + \theta^{-1}$ has a multiplicative order at least $U(r-2, p-1)$ and this order is a divisor of $q^{(r-1)/2}-1$*

*Proof.* A fact that a multiplicative order of element $\theta+\theta^{-1} = \theta^{-1}(\theta^2+1)$ is at least $U(r-2,p-1)$ follows from the theorem 1, (a). Since

$$(\theta+\theta^{-1})^{q^{(r-1)/2}-1} = (\theta^{q^{(r-1)/2}} + \theta^{-q^{(r-1)/2}})(\theta+\theta^{-1})^{-1} = (\theta^{-1}+\theta)(\theta+\theta^{-1})^{-1} = 1$$

an order of element $\theta+\theta^{-1}$ is a divisor of $q^{(r-1)/2}-1$. □

**Corollary 3.** *If $a \neq 1, -1$ then element*

$$z = \begin{cases} (\theta+\theta^{-1})^2(a\theta+1)(\theta+a)^{-1} & \text{if } \rho_2(q^{(r-1)/2}-1) = 2 \\ (\theta+\theta^{-1})[(a\theta+1)(\theta+a)^{-1}]^2 & \text{if } \rho_2(q^{(r-1)/2}+1) = 2 \end{cases}$$

*has a multiplicative order at least $[U(r-2, p-1)\, U((r-3)/2, p-1)]/2$*

*Proof.* According to corollary 2 element $\theta+\theta^{-1}$ has order that divides $q^{(r-1)/2}-1$ and generates subgroup of order at least $U(r-2, p-1)$, element $(\theta+\theta^{-1})^2$ has order that divides $(q^{(r-1)/2}-1)/2$ and generates subgroup of order at least $U(r-2, p-1)/2$. According to proof of the theorem 1, (b) element $(a\theta+1)(\theta+a)^{-1}$ has order that divides $q^{(r-1)/2}+1$ and generates subgroup of order at least $U((r-3)/2, p-1)$, element $[(a\theta+1)(\theta+a)^{-1}]^2$ has order that divides $(q^{(r-1)/2}+1)/2$ and generates subgroup of order at least $U((r-3)/2, p-1)/2$.

If $\rho_2(q^{(r-1)/2}-1)=2$ then orders of elements $(\theta+\theta^{-1})^2$ and $(a\theta+1)(\theta+a)^{-1}$ are coprime.

If $\rho_2(q^{(r-1)/2}+1)=2$ then orders of elements $\theta+\theta^{-1}$ and $[(a\theta+1)(\theta+a)^{-1}]^2$ are coprime.

Both in the first and in the second case an order of element $z$ is a product of orders of its factors. $\square$

## 3  Explicit lower bounds on multiplicative orders

We use some known estimates from [3,7] to derive explicit lower bounds on the multiplicative order of the elements $\theta^e(\theta^f+a)$ and $z$.

According to [3, corollary 1.3(Glaisher)] the following equality is true :

$$U(n,d-1)=Q(n,d) \qquad (13)$$

We consider two different cases.

Case 1) $r-3 \geq 2p^2$, that is $r$ is big comparatively to $p$.

In this case the following corollary holds. Note that $2.5 \approx \pi\sqrt{2/3}$.

**Corollary 4.** *Let $r-3\geq 2p^2$, $e$ be any integer, $f$ be any integer coprime with $r$, $a\in F_q$ be any non-zero element. Then*

*(a) element $\theta^e(\theta^f+a)$ in the extension field $F_q(\theta)$ has a multiplicative order larger than*

$$\left(\frac{p(p-1)}{160(r-2)}\right)^{\sqrt{p}} \exp\left(2.5\sqrt{(1-\frac{1}{p})(r-2)}\right)$$

*(b) if $a\neq 1,-1$ then element $\theta^e(\theta^f+a)$ has a multiplicative order larger than*

$$\frac{1}{2}\left(\frac{p(p-1)}{80(r-3)}\right)^{2\sqrt{p}} \exp\left(2.5\sqrt{2}\sqrt{(1-\frac{1}{p})(r-3)}\right)$$

*(c) if $a\neq 1,-1$ then element*

$$z = \begin{cases} (\theta+\theta^{-1})^2(a\theta+1)(\theta+a)^{-1} & \text{if } \rho_2(q^{(r-1)/2}-1)=2 \\ (\theta+\theta^{-1})[(a\theta+1)(\theta+a)^{-1}]^2 & \text{if } \rho_2(q^{(r-1)/2}+1)=2 \end{cases}$$

*has a multiplicative order larger than*

$$\frac{1}{2}\left(\frac{p(p-1)}{80(r-3)}\right)^{\sqrt{p}} \exp\left(2.5(1+\frac{\sqrt{2}}{2})\sqrt{(1-\frac{1}{p})(r-3)}\right)$$

*Proof.*

(a) According to [7, theorem 5.1] the following inequality holds for $n\geq d^2$

$$Q(n,d) > \left(\frac{d(d-1)}{160n}\right)^{\sqrt{n}} \exp\left(2.5\sqrt{(1-\frac{1}{d})n}\right) \qquad (14)$$

According to theorem 1, equality (13) and inequality (14) we have for $r-2 \geq p^2$:

$$ord(\theta^e(\theta^f+a)) \geq U(r-2, p-1) = Q(r-2, p) > \left(\frac{p(p-1)}{160(r-2)}\right)^{\sqrt{p}} \exp\left(2.5\sqrt{(1-\frac{1}{p})(r-2)}\right)$$

(b) Analogous to proof of (a) using [7,theorem 5.1], theorem, equality (13) and inequality (14). Note that if $r-3 \geq 2p^2$ then $r-2 \geq p^2$.

(c) Analogous to proof of (a) using [2,theorem 5.1], theorem, equality (13) and inequality (14) □

Case 2) $r-2<p$, that is $r$ is the same magnitude as $p$ or small comparatively to $p$

In this case the following corollary holds.

**Corollary 5.** *Let $r-2<p$, e be any integer, f be any integer coprime with r, $a \in F_q$ be any non-zero element. Then*

*(a) element $\theta^e(\theta^f+a)$ in the extension field $F_q(\theta)$ has a multiplicative order larger than*

$$\frac{\exp(2.5\sqrt{r-2})}{13(r-2)}$$

*(b) if $a \neq 1,-1$ then element $\theta^e(\theta^f+a)$ has a multiplicative order larger than $\frac{2\exp(2.5\sqrt{2}\sqrt{r-3})}{169(r-3)^2}$*

*(c) if $a \neq 1,-1$ then element*

$$z = \begin{cases} (\theta+\theta^{-1})^2(a\theta+1)(\theta+a)^{-1} & \text{if } \rho_2(q^{(r-1)/2}-1) = 2 \\ (\theta+\theta^{-1})[(a\theta+1)(\theta+a)^{-1}]^2 & \text{if } \rho_2(q^{(r-1)/2}+1) = 2 \end{cases}$$

*has a multiplicative order larger than*

$$\frac{\exp\left(2.5(1+\frac{\sqrt{2}}{2})\sqrt{r-3}\right)}{169(r-2)(r-3)}$$

*Proof.*

(a) If $n<p$ then $U(n,p-1)=U(n)$. According to [7,theorem 4.2] the following inequality holds for $n<d$

$$U(n) > \frac{\exp(2.5\sqrt{n})}{13n} \qquad (15)$$

According to theorem 1, equality 13 and inequality (15) we have for $r-2<p$:

$$ord(\theta^e(\theta^f+a)) \geq U(r-2, p-1) = U(r-2) > \frac{\exp(2.5\sqrt{r-2})}{13(r-2)}$$

(b) Analogous to proof of (a) using [7,theorem 4.2], theorem, equality (13) and inequality (15). Note that if $r-2<p$ then $(r-3)/2<p$.

(c) Analogous to proof of (a) using [7,theorem 4.2], theorem, equality (13) and inequality (15) □

**Remark 2.** *Element $\theta^e(\theta^f + a)$ asymptotically has a multiplicative order larger than*

$$\exp(2.5\sqrt{r}) = 12.18^{\sqrt{r}}.$$

*If $a \neq 1, -1$ then element $\theta^e(\theta^f + a)$ asymptotically has a multiplicative order larger than*

$$\exp(2.5\sqrt{2}\sqrt{r}) = 33.95^{\sqrt{r}}.$$

*If $a \neq 1, -1$ then element*

$$z = \begin{cases} (\theta + \theta^{-1})^2(a\theta + 1)(\theta + a)^{-1} & \text{if } \rho_2(q^{(r-1)/2} - 1) = 2 \\ (\theta + \theta^{-1})[(a\theta + 1)(\theta + a)^{-1}]^2 & \text{if } \rho_2(q^{(r-1)/2} + 1) = 2 \end{cases}$$

*asymptotically has a multiplicative order larger than*

$$\exp\left(2.5(1 + \frac{\sqrt{2}}{2})\sqrt{r}\right) = 70.1^{\sqrt{r}}$$

## 4    Examples

Let us denote lower bounds for orders of elements $\theta^e(\theta^f + 1)$, $\theta^e(\theta^f + a)$,

$$z = \begin{cases} (\theta + \theta^{-1})^2(a\theta + 1)(\theta + a)^{-1} & \text{if } \rho_2(q^{(r-1)/2} - 1) = 2 \\ (\theta + \theta^{-1})[(a\theta + 1)(\theta + a)^{-1}]^2 & \text{if } \rho_2(q^{(r-1)/2} + 1) = 2 \end{cases}$$

by $z_1$, $z_2$, $z_3$ respectively.

Logarithms of $|F^*_{q^{r-1}}|$ and of $z_1$, $z_2$, $z_3$ in examples 1-4 are given in the table.

**Example 1**

$q=p=5$, $r=257$ – prime number, element 5 is primitive modulo 257 and $F_{q^{r-1}} = F_{5^{256}}$. Since $r-3 \geq 2p^2$ we have case 1 in this example.

**Example 2**

$q=p=3$, $r=401$ – prime number, element 3 is primitive modulo 401 and $F_{q^{r-1}} = F_{3^{400}}$. Since $r-3 \geq 2p^2$ we have case 1 in this example.

**Example 3**

$q=p=11$, $r=1009$ – prime number, element $q=11$ is primitive modulo $r=1009$ and $F_{q^{r-1}} = F_{11^{1008}}$. Since $r-3 \geq 2p^2$ we have case 1 in this example.

**Example 4**

$q=p=107$, $r=97$ – prime number, element $q=107$ is primitive modulo $r=97$ and $F_{q^{r-1}} = F_{107^{96}}$. Since $r-2<p$ we have case 2 in this example.

| | Q | r | $\log_2 |F^*_{q^{r-1}}|$ | $\log_2 z_1$ | $\log_2 z_2$ | $\log_2 z_3$ |
|---|---|---|---|---|---|---|
| 1 | 5 | 257 | 594.41 | 26.93 | 27.03 | 64.43 |
| 2 | 3 | 401 | 634 | 35.65 | 39.22 | 77.86 |
| 3 | 11 | 1009 | 3487.1 | 74.24 | 90.13 | 153.64 |
| 4 | 107 | 97 | 647.18 | 24.71 | 28.71 | 38.89 |

Department of Computer Technologies,
National University Lviv Politechnika,
Bandery Str.,12, 79013, Lviv, Ukraine
e-mail: popovych@polynet.lviv.ua